\documentclass[a4paper,12pt]{article}
\usepackage{amsfonts,amsmath,amssymb,amscd,latexsym}
\usepackage{srcltx}

\linethickness{0.4pt}

\newcommand{\bi}{\bibitem}
\newcommand{\nb}{\newblock}

\newcommand{\be}[1]{\begin{equation}\label{#1}}
\newcommand{\ee}{\end{equation}}

\newcommand{\la}{\langle\,}
\newcommand{\ra}{\,\rangle}
\newcommand{\ve}{\varepsilon}
\newcommand{\prf}{{\bf Proof.}\ }

\newcommand{\pp}{{\mathcal P}}
\newcommand{\rr}{{\mathcal R}}
\newcommand{\dd}{{\mathcal D}}

\newcommand{\topp}{\mathop{\mbox{\bf top}}}
\newcommand{\bott}{\mathop{\mbox{\bf bot}}}

\newtheorem{thm}{\quad Theorem}
\newtheorem{lm}{\quad Lemma}

\begin{document}

\title{On diagram groups over Fibonacci-like semigroup presentations and their generalizations}

\author{\vspace{2ex}
V. S. Guba\thanks{This work is partially supported by the Russian Foundation
for Basic Research, project no. 19-01-00591 A.}\\
Vologda State University,\\
15 Lenin Street,\\
Vologda\\
Russia\\
160600\\
E-mail: guba{@}uni-vologda.ac.ru}
\date{}

\maketitle

\begin{abstract}

We answer the question by Matt Brin on the structure of diagram groups over semigroup
presentation ${\mathcal P}=\langle a,b,c\mid a=bc,b=ca,c=ab\rangle$. In the talk on Oberwolfach workshop,
Brin conjectured that the diagram group over $\mathcal P$ with base $a$ is isomorphic to the generalized
Thompson's group $F_9$. We confirm this conjecture and consider some generalizations of this fact.

\end{abstract}

\section{Introduction}
\label{backgr}

In this background Section we recall the concept of diagram groups and introduce some terminology. The contents of the present
Section is essentially known. Some defininions and examples from here repeat the ones from \cite{Gu04}. Detailed information about
diagram groups can be found in \cite{GbS}.

First of all, let us recall the concept of a semigroup diagram and
introduce some notation. To do this, we consider the following example.
Let $\pp=\la a,b\mid aba=b,bab=a\ra$ be the semigroup presentation. (In the next Section we will work with it.)

It is easy to see by the following algebraic calculation
$$
a^5=a(bab)a(bab)a=(aba)(bab)(aba)=bab=a
$$
that the words $a^5$ and $a$ are equal modulo $\pp$. The same can be seen from the following picture

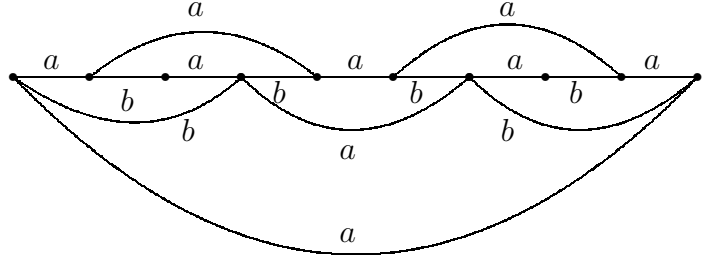
\begin{center}
\begin{picture}(90.00,37.00)
\put(00.00,23.00){\circle*{1.00}}
\put(10.00,23.00){\circle*{1.00}}
\put(20.00,23.00){\circle*{1.00}}
\put(30.00,23.00){\circle*{1.00}}
\put(30.00,23.00){\circle*{1.00}}
\put(40.00,23.00){\circle*{1.00}}
\put(50.00,23.00){\circle*{1.00}}
\put(60.00,23.00){\circle*{1.00}}
\put(60.00,23.00){\circle*{1.00}}
\put(70.00,23.00){\circle*{1.00}}
\put(80.00,23.00){\circle*{1.00}}
\put(90.00,23.00){\circle*{1.00}}
\put(00.00,23.00){\line(1,0){90.00}}
\bezier{152}(10.00,23.00)(25.00,35.00)(40.00,23.00)
\bezier{240}(50.00,23.00)(80.00,23.00)(50.00,23.00)
\bezier{164}(50.00,23.00)(65.00,37.00)(80.00,23.00)
\bezier{240}(0.00,23.00)(30.00,23.00)(0.00,23.00)
\bezier{156}(0.00,23.00)(17.00,11.00)(30.00,23.00)
\bezier{164}(30.00,23.00)(44.00,9.00)(60.00,23.00)
\bezier{164}(60.00,23.00)(74.00,9.00)(90.00,23.00)
\put(5.00,25.00){\makebox(0,0)[cc]{$a$}}
\put(24.00,32.00){\makebox(0,0)[cc]{$a$}}
\put(45.00,25.00){\makebox(0,0)[cc]{$a$}}
\put(65.00,32.00){\makebox(0,0)[cc]{$a$}}
\put(84.00,25.00){\makebox(0,0)[cc]{$a$}}
\put(23.00,16.00){\makebox(0,0)[cc]{$b$}}
\put(44.00,13.00){\makebox(0,0)[cc]{$a$}}
\put(65.00,16.00){\makebox(0,0)[cc]{$b$}}
\put(15.00,20.00){\makebox(0,0)[cc]{$b$}}
\put(24.00,25.00){\makebox(0,0)[cc]{$a$}}
\put(35.00,21.00){\makebox(0,0)[cc]{$b$}}
\put(53.00,21.00){\makebox(0,0)[cc]{$b$}}
\put(66.00,25.00){\makebox(0,0)[cc]{$a$}}
\put(74.00,21.00){\makebox(0,0)[cc]{$b$}}
\bezier{520}(0.00,23.00)(45.00,-24.00)(90.00,23.00)
\put(44.00,2.00){\makebox(0,0)[cc]{$a$}}
\end{picture}
\end{center}

This is a {\em diagram\/} $\Delta$ over the semigroup presentation $\pp$.
It is a plane graph with $10$ vertices, $15$ (geometric) edges and $6$ faces
or {\em cells\/}. Each cell corresponds to an elementary transformation of a
word, that is, a transformation of the form
$p\cdot u\cdot q\to p\cdot v\cdot q$, where $p$, $q$ are words (possibly,
empty), $u=v$ or $v=u$ belongs to the set of defining relations. The
diagram $\Delta$ has the leftmost vertex denoted by $\iota(\Delta)$ and
the rightmost vertex denoted by $\tau(\Delta)$. It also has the {\em top
path\/} $\topp(\Delta)$ and the {\em bottom path\/} $\bott(\Delta)$ from
$\iota(\Delta)$ to $\tau(\Delta)$. Each cell $\pi$ of a diagram can be
regarded as a diagram itself. The above functions $\iota$, $\tau$, $\topp$,
$\bott$ can be applied to $\pi$ as well. We do not distinguish isotopic
diagrams.

We say that $\Delta$ is a $(w_1,w_2)$-diagram whenever the label of its
top path is $w_1$ and the label of its bottom path is $w_2$. In our example,
we deal with an $(a^5,a)$-diagram. If we have two diagrams such that the
bottom path of the first of them has the same label as the top path of the
second, then we can naturally {\em concatenate\/} these diagrams by
identifying the bottom path of the first diagram with the top path of the
second diagram. The result of the concatenation of a $(w_1,w_2)$-diagram and
a $(w_2,w_3)$-diagram obviously is a $(w_1,w_3)$-diagram. We use the sign
$\circ$ for the operation of concatenation. For any diagram $\Delta$ over
$\pp$ one can consider its {\em mirror image\/} $\Delta^{-1}$. A diagram may
have {\em dipoles\/}, that is, subdiagrams of the form $\pi\circ\pi^{-1}$,
where $\pi$ is a single cell. To {\em cancel\/} (or {\em reduce\/}) the
dipole means to remove the common boundary of $\pi$ and $\pi^{-1}$ identifying $\topp(\pi)$ with $\bott(\pi^{-1})$.
In any diagram, we can cancel all its dipoles, step by step. The result does not depend on the order of
cancellations. A diagram is {\em irreducible\/} whenever it has no dipoles.
The operation of cancelling dipoles has an inverse operation called the
{\em insertion\/} of a dipole. These operations induce an equivalence
relation on the set of diagrams (two diagrams are {\em equivalent\/} whenever
one can go from one of them to the other by a finite sequence of
cancelling/inserting dipoles). Each equivalence class contains exactly one
irreducible diagram.

For any nonempty word $w$, the set of all $(w,w)$-diagrams forms a monoid
with the identity element $\ve(w)$ (the diagram with no cells). The
operation $\circ$ naturally induces some operation on the set of equivalence
classes of diagrams. This operation is called a {\em product\/} and
equivalent diagrams are called {\em equal\/}. (The sign $\equiv$ will be
used to denote that two diagrams are isotopic.) So the set of all
equivalence classes of $(w,w)$-diagrams forms a group that is called the
{\em diagram group\/} over $\pp$ with {\em base\/} $w$. We denote this
group by ${\mathcal D}(\pp,w)$. We can think of this group as of the set
of all irreducible $(w,w)$-diagrams. The group operation is the
concatenation with cancelling all dipoles in the result. An inverse
element of a diagram is its mirror image. We also need one more natural
operation on the set of diagrams. By the {\em sum\/} of two diagrams we
mean the diagram obtained by identifying the rightmost vertex of the
first summand with the leftmost vertex of the second summand. This
operation is also associative. The sum of diagrams $\Delta_1$,
$\Delta_2$ is denoted by $\Delta_1+\Delta_2$.
\vspace{2ex}

Now let us recall some information about generalized Thompson's groups $F_r$.
This family was introduced by K. S. Brown in \cite{Bro}. Additional facts about these groups
can be found in \cite{BCS,Stein}.

The family of generalized Thompson's groups can be defined as follows.
The group $F_r$ is the group of all piecewise linear self homeomorphisms
of the unit interval $[0,1]$ that are orientation preserving (that is,
send $0$ to zero and $1$ to $1$) with all slopes integer powers of $r$
and such that their singularities (breakpoints of the derivative) belong
to $\mathbb Z[\,\frac1r\,]$. The group $F_r$ admits a presentation given by
\be{presfp}
\la x_0,x_1,x_2,\ldots\mid x_jx_i=x_ix_{j+r-1}\ (i<j)\ra.
\ee

\noindent
This presentation is infinite, but a close examination shows that the
group is actually finitely generated, since $x_0$, $x_1$, \dots, $x_{p-1}$
are sufficient to generate it. In fact, the group is finitely presented; see \cite{Bro}.
The finite presentation is awkward, and it is not used much. The symmetric
and simple nature of the infinite presentation makes it much more adequate
for almost all purposes.

One way in which the infinite presentation is very useful is in the construction of the normal forms.
A word given in the generators $x_i$ and their inverses, can have its generators moved around
according to the relators, and the result is the following well-known statement:

\begin{thm}
\label{stnf}
An element in $F_r$ always admits an expression of the form
$$
x_{i_1}x_{i_2}\cdots x_{i_m}x_{j_n}^{-1}\cdots x_{j_2}^{-1}x_{j_1}^{-1},
$$
where
$$
i_1\le i_2\le\cdots\le i_m,\ j_1\le j_2\le\cdots\le j_n.
$$
\end{thm}

In general, this expression is not unique, but for every element there is
a unique word of this type which satisfies certain technical condition.
This unique word is called the {\em standard normal form\/} for the
element of $F_r$.

The case $r=2$ corresponds to famous R. Thompson's group $F=F_2$.

It is known \cite{GbS} that groups $F_r$ are diagram groups over the
semigroup presentation $\pp_r=\la x\mid x=x^r\ra$ with base $x$ (note that
for any base $x^k$, where $k\ge1$, we get an isomorphic group).

Now let us compare the diagram representation of $F$ with the
representation of its elements by piecewise-linear homeomorphisms of
the closed unit interval $[0,1]$. Let $\Delta$ be an $(x^p,x^q)$-diagram
over $\pp$. We will show how to assign to it a piecewise-linear function
from $[0,p]$ onto $[0,q]$. Each positive edge of $\Delta$ is homeomorphic
to the unit interval $[0,1]$. So we assign a coordinate to each point of
this edge (the leftmost end of an edge has coordinate $0$, the rightmost one
has coordinate $1$). Let $\pi$ be an $(x,x^r)$-cell of $\Delta$. Let us map
$\topp(\pi)$ onto $\bott(\pi)$ linearly, that it, the point on the edge
$\topp(\pi)$ with coordinate $t\in[0,1]$ is taken to the point on
$\bott(\pi)$ with coordinate $rt$ (the bottom path of $\pi$ has length
$r$ so it is naturally homeomorphic to $[0,r]$). The same thing can be
done for an $(x^r,x)$-cell of $\Delta$. Thus for any cell $\pi$ of $\Delta$
we have a natural mapping $T_\pi$ from $\topp(\pi)$ onto $\bott(\pi)$ (we
call it a {\em transition map\/}). Now let $t$ be any number in $[0,p]$. We
consider the point $o$ on $\topp(\Delta)$ that has coordinate $t$. If $o$ is
not a point of $\bott(\Delta)$, then it is an internal point on the top path
of some cell. Thus we can apply the corresponding transition map to $o$. We
repeat this operation until we get a point $o'$ on the path $\bott(\Delta)$.
The coordinate of this point is a number in $[0,q]$. Hence we have a function
$f_\Delta\colon[0,p]\to[0,q]$ induced by $\Delta$. It is easy to see this
will be a piecewise-linear function. When we concatenate diagrams, this
corresponds to the composition of the PL functions induced by these diagrams.
For groups $F_r$, which are the diagram group $\dd(\pp_r,x)$, we have the
homomorphism from it to $PLF[0,1]$. It is known this is an monomorphism.
\vspace{1ex}

The following elementary fact was essentially used several times in \cite{GuSa99,Gu00}
and some other papers.

\begin{lm}
\label{longpath}
Let $\pp=\la X\mid\rr\ra$ be a semigroup presentation. Suppose that all
defining relations of $\pp$ have the form $a=A$, where $a\in X$ and $A$ is
a word of length at least $2$. Also assume that all letters in the
left-hand sides of the defining relations are different. Then any
irreducible diagram $\Delta$ over $\pp$ is the concatenation of the form
$\Delta_1\circ\Delta_2^{-1}$, where the top path of each cell of both
$\Delta_1$, $\Delta_2$ has length $1$. The longest positive path in $\Delta$
from $\iota(\Delta)$ to $\tau(\Delta)$ coincides with the bottom path of
$\Delta_1$ and the top path of $\Delta_2^{-1}$.
\end{lm}

Note that $\la x\mid x=x^r\ra$ obviously satisfies the conditions of the
Lemma. The same concerns the presentation $\la a,b\mid a=bab, b=aba\ra$,
which was considered in the beginning of this Section. Let us recall the
idea of the proof. Let $p$ be the longest positive path in $\Delta$ from
$\iota(\Delta)$ to $\tau(\Delta)$. It cuts $\Delta$ into two parts. It
suffices to prove that all cells in the ``upper" part correspond to the
defining relations of the form $a=A$, where $a$ is a letter, and none of
them corresponds to $A=a$. Assume the contrary. Suppose that there is a
cell $\pi$ in the upper part of $\Delta$ with the top label $A$ and the
bottom label $a$. The bottom path of $\pi$ cannot be a subpath in $p$ since
$p$ is chosen the longest. So the bottom edge of $\pi$ belongs to the top
path of some cell $\pi'$. The diagram $\Delta$ has no dipoles. All letters
in the left-hand sides of the defining relations are different. So the top
path of $\pi'$ cannot have length $1$. This means that we have found a new
cell in the upper part of $\Delta$ that also corresponds to the defining
relation of the form $A=a$. Applying the same argument to $\pi'$, we get a
process that never terminates. This is impossible since the cells that appear
during the process cannot repeat. This completes the proof.

\section{Main Results}
\label{flp}

Let $a_1$, $a_2$, ... , $a_n$ be a finite alphabet. By definition, $a_{n+1}=a_1$, $a_{n+2}=a_2$. Consider the
following semigroup presentation
\be{fs}
{\mathcal P}_n=\langle a_1,\ldots,a_n\mid a_i=a_{i+1}a_{i+2}\ (1\le i\le n)\rangle.
\ee
The semigroup presented by ${\mathcal P}_n$ is called {\em Fibonacci semigroup\/}. One can ask what are the
diagram groups $G_n=\mathcal D(\mathcal P_n,a_1)$. The case $n=1$ is trivial, it gives the diagram group over
$\langle x\mid x=xx\rangle$ so it is Thompson's group $F$. For $n=2$ one has the presentation $\langle a,b\mid
a=ba,b=ab\rangle$. It was shown in \cite{GoSa17} that $G_2$ (the so called Jones' subgroup) is isomorphic to
$F_3$.

In his talk on an Oberwolfach worksop, Matt Brin asked about the group $G_3$, see \cite[Question 73]{BBN}. He conjectured
that this diagram group is isomorphic to $F_9$. Notice that $\mathcal P_3$ can be written as $\langle a,b,c\mid
a=bc,b=ca,c=ab\rangle$. This presentation is not {\em complete\/}. This means that for the Thue system $ab\to
c$, $bc\to a$, $ca\to b$ there are no unique normal forms. For complete semigroup presentations, there exists a
technique of their calculation from \cite{GbS}. Sometimes it is possible to consider a completion, but here it
has a complicated form. Indeed, the semigroup given by $\mathcal P_3$ is the quaternion group $Q_8$. So this way
of description looks very unclear.

Here we present a purely geometric way to find the diagram group. First of all, let us mention that one can
avoid generator $c$ replacing it by $ab$. In generators $a$, $b$ the presentaion becomes $\langle a,b\mid
a=bab,b=aba\rangle$. It was considered as an example in the beginning of the Introduction. The semigroup
given by it is the same as above.

There is a fact from \cite[Section 4]{GuSa05} that ordinary Tietze transformations of semigroup presentations
lead to the same diagram groups. (This can also be shown directly.) So we have one more generalization of the
class of semigroup presentations under consideration.

Let $a_1$, $a_2$, ... , $a_n$ be a finite alphabet as above and let $r\ge2$ be an integer. For any $j$ from $1$
to $r$ we set $a_{n+j}=a_j$. Now for every $i$ from $1$ to $n$ we consider a relation of the form
$a_i=a_{i+1}\ldots a_{i+r}$. By $\mathcal P_{nr}$ we denote a semigroup presentation given by these relations:
\be{pnr}
{\mathcal P}_{nr}=\langle a_1,\ldots,a_n\mid a_i=a_{i+1}\ldots a_{i+r}\ (1\le i\le n)\rangle.
\ee

This class of presentations was introduced by Johnson in \cite{Jo74} in order to generalize the concept of a
Fibonacci group. Since $\mathcal P_{nr}$ is also a semigroup presentation, one can introduce the corresponding
semigroups as well. For $r=2$ we have the above Fibonacci-like presentations. Now we can consider diagram groups
$G_{nr}$ defined as $\mathcal D(\mathcal P_{nr},a_1)$. The group we are interested in is $G_{32}\cong G_{23}$.
We confirm Brin's conjecture about it.

\begin{thm}
\label{f9}
The diagram group with base $a$ over semigroup presentation $\langle a,b,c\mid a=bc,b=ca,c=ab\rangle$ is
isomorphic to generalized Thompson's group $F_9$.
\end{thm}

\prf
We consider this group as a diagram group over $\mathcal P_{23}=\langle a,b\mid a=bab,b=aba\rangle$. It is known
that the groups $F_r$ have no proper non-Abelian homomorphic images. So it suffices to construct a homomorphism
from $F_9$ to the diagram group $G=G_{23}$ showing it is surjective. Therefore, this will give us an
isomorphism.

The group $F_9$ will be considered as the diagram group over $\langle x\mid x=x^9\rangle$ with base $x$. A
diagram over this presentation is a plane graph composed from cycles of even length. By induction on the number
of cells it is easy to show that the graph is bipartite. So we can give colours to its vertices. Let the initial
vertex of a diagram $\Delta$ gets the colour 1. Then the other vertices get their colours uniquely.

Now we relabel the diagram: if a positive edge goes from a vertex of colour 1 to the vertex of color 2, then we
give it label $a$. Otherwise it has label $b$. As a result, we get a diagram denoted by $\Delta'$. Each
cell $x=x^9$ becomes a cell of one of the two forms: $a=a(ba)^4$ or $b=b(ab)^4$. The same for inverse cells.

We have the following derivation over $\mathcal P_{23}$: $a=bab=(aba)(bab)(aba)=a(ba)^4$, and similarly for the
other equality. Semigroup diagrams for these equalities consist of $4$ cells. They will be called {\em basic\/}. We fill the cells
of the above form by basic diagrams. This gives us the diagram $\Delta''$ over $\mathcal P_{23}$.

The rule $\Delta\mapsto\Delta''$ induces a homomorphism of groupoids of diagrams. (Notice that cancelling a
dipole in a diagram $\Delta$ over $x=x^9$ leads to cancelling $4$ dipoles in $\Delta''$ so the mapping we have
defined preserves equivalence of diagrams.) In particular, we have a homomorphism from $F_9$ as the diagram group
over $x=x^9$ with base $x$ to $G$ as the diagram over $\mathcal P_{23}$ with base $a$.

Now let $\Psi$ be a reduced $(a,a)$-diagram over $\mathcal P_{23}$. We would like to find a preimage of it in
$F_9$. According to Lemma~\ref{longpath}, we decompose $\Psi$ as $\Psi_1\circ\Psi_2^{-1}$ where $\Psi_1$,
$\Psi_2$ are positive diagrams. It holds that $\bott{\Psi_1}=\topp{\Psi_2^{-1}=p}$, where $p$ is the longest
positive path in $\Psi$ from $\iota(\Psi)$ to $\tau(\Psi)$.

Now we will change $\Psi=\Psi_1\circ\Psi_2^{-1}$ and the path $p$ step by step inserting some dipoles. The
current situation will always have the same notation. Suppose that the first edge of $p$ has label $b$. In this
case we replace the subdiagram $\varepsilon(b)$ that consists of one edge by a dipole of the two cells
$(b=aba)\circ(aba=b)$. The new longest path in the diagram we obtain will be still denoted by $p$.

Now look and the subwords of the form $aa$ or $bb$ of the label of $p$. Choose the leftmost of them. If it is
$aa$ then we replace the second edge labelled by $a$ by the dipole $(a=bab)\circ(bab=b)$. If it is $bb$ then
we also replace the second edge of it by the dipole $(b=aba)\circ(aba=b)$.

After a finite number of steps, the label of the longest path $p$ becomes $abab\ldots$\ . The last letter in it
will have label $a$. This follows from parity arguments and the fact that the terminal vertex of $\Psi$ has
colour $2$. Now we have $\Psi=\Psi_1\circ\Psi_2^{-1}$ where $\Psi_1$, $\Psi_2$ are positive $(a,a(ba)^m)$-diagrams
for some $m$. It suffices to show that each diagram with this property belongs to the image of our mapping
$\Delta\mapsto\Delta''$. This means that every positive $(a,a(ba)^m)$-diagram over $\mathcal P_{23}$ can be
composed from basic diagrams. Also we claim a symmetric statement: every positive $(b,b(ab)^m)$-diagram over $\mathcal P_{23}$
can be composed from basic diagrams.

Let $\Phi$ be one of these diagrams. We proceed by induction on the number of cells in it. If there are no cells
($m=0$) then we can nothing to prove. Otherwise let us define the {\em depth\/} of an edge in the diagram. The
top edge will have depth $0$ by definition. All other edges belong to the bottom path of a cell $\pi$. If its top edge has depth $d$,
then we assign depth $d+1$ to our edge. The only important thing for us is whether $d$ is even or odd. So we talk about even and odd edges.

Now we remark the following.

1) Let $e_1$, ... , $e_s$ be all edges coming out of a vertex, read from top to bottom. Then labels of them always change from
$a$ to $b$ and vice versa, and the same for parity of their depth. The same for edges that come into a vertex.

2) If two consecutive edges have the same label, then they have different partity. Otherwise, if the labels are $ab$ or $ba$,
the parity is the same.

The first part is clear. As for the second one, let us consider only one case of the edges labelled by $ab$. Let $e$ be the highest
edge that ends at $v$ (the vertex between $a$ and $b$) and let $f$ be the highest edge that starts at $v$. It is easy to see that $ef$ is a
part of the bottom path of a cell. Therefore, $e$ and $f$ have different labels and the same depth. Now everything follows from 1).
The cases $ba$, $aa$, $bb$ are similar.

Now we look again at the path $p$ (the bottom of $\Phi$). Its first label is $a$, so the first edge is even. Therefore,
all edges of $p$ are even according to 2) since $p$ has label $abab...a$. If $m > 0$ then $\Phi$ has a top cell $a=bab$ with the bottom
path $e_1e_2e_3$. Deleting the top cell gives us a sum of 3 diagrams: $(e_1,p_1)+(e_2,p_2)+(e_3,p_3)$, where $p=p_1p_2p_3$. Each edge
of $p_i$ has odd parity in the $i$-th sumand. Therefore, $e_i$ does not belong to $p_i$. So there exists a top cell in each of the
summands. Together with the cell $a=bab$ we have deleted, they form a basic diagram.

Removing three cells with top edges $e_i$ ($i=1,2,3)$, we get a sum of $9$ positive diagrams. Now all edges of $p$ have even depth so the
inductive assumption can be applied to these summands. This completes the proof.
\vspace{2ex}

So this answers Brin's question, and now we look at some generalizations. The next Fibonacci-like presentation
in the series is $\mathcal P_{42}$. Its relations are $a=bc$, $b=cd$, $c=da$, $d=ab$. Applying Tietze
transformations, we rewrite the presentation as $\langle a,b\mid a=baba,b=abaab\rangle$, where $d\to ab$, $c\to da\to
aba$. In the second relation $b=abaab$ we replace its third occurrence of $a$ to the right-hand side by
$(ba)^2$. This gives us a Tietze-equivalent presentation $\mathcal P=\langle a,b\mid a=(ba)^2,b=(ab)^4$. The diagram
group over $\mathcal P$ with base $a$ is the same as the one over $\mathcal P_{42}$ according to general facts from
\cite{GuSa05}.

\begin{thm}
\label{f11}
The diagram group with base $a$ over semigroup presentation $\langle a,b,c\mid a=bc,b=cd,c=da,d=ab\rangle$ is
isomorphic to generalized Thompson's group $F_{11}$.
\end{thm}

\prf The idea of the proof is similar to the one for Theorem~\ref{f9}. We will work with presentation $\mathcal P=\langle a,b\mid a=(ba)^2,b=(ab)^4$
instead of $\mathcal P_{42}$. Our aim is to construct a homomorphism from $F_{11}$ to $G=\mathcal D(\mathcal
P,a)$. Notice that we have no longer a symmetry between $a$ and $b$. The group $F_{11}$ will be the diagram
group with base $x$ over $x=x^{11}$, as usual. Any diagram $\Delta$ over it is still a bipartite graph since
$11$ is odd. So each vertex gets a colour $1$ or $2$ and each edge will have a label $a$ or $b$ by the same
rules as above. This new diagram over $a=a(ba)^5$, $b=b(ab)^5$ will be denoted by $\Delta'$.

Both relations can be derived from $\mathcal P$. Indeed, $a=baba$, and then we replace the first occurrence of
$b$ to the right-hand side by $(ab)^4$. Thus we have a diagram of two cells over $\mathcal P$ for $a=a(ba)^5$. As for the second
equality, we take $b=ababab$ and replace the first $a$ by $(ba)^2$. This gives a two-cell diagram over $\mathcal P$ for
$b=b(ab)^5$. These two diagrams over $\mathcal P$ will be called basic. Replacing the cells of $\Delta'$ by
basic diagrams lead to the diagram $\Delta''$. In a standard way, the mapping $\Delta\mapsto\Delta''$ induces
the homomorphism of the groupoids of diagrams, and therefore we have a group homomorphism from $F_{11}$ to $G$.
Our aim is to establish its surjectivity.

Now let $\Psi$ be a reduced diagram over $\mathcal P$. As in the proof of the previous theorem, we let
$\Psi=\Psi_1\circ\Psi_2^{-1}$ where $p$ is the common part of the two pieces. We are going to insert certain
dipoles to $\Psi$ in such a way that the label of $p$ will have the form $abab\ldots$\ .

Suppose that the label of $p$ starts with $b$. Then we insert a dipole of the form $(b=(ab)^4)\circ((ab)^4=b)$
instead of the first edge of $p$. The new path is still denoted by $p$. If its label has an occurrence of $aa$
or $bb$ then we take the leftmost of them. In case it is $aa$, we replace the second edge by the dipole
$(a=(ba)^2)\circ((ba)^2=a)$. In case it is $bb$, the second edge is replaced by a dipole from the beginning of this
paragraph. So in a finite number of steps, we get a decomposition into a product of two diagrams, positive
and negative. It suffices to take a positive $(a,abab\ldots)$-diagram $\Phi$ showing that it is in the image
of the mapping $\Delta\mapsto\Delta''$.

Now we are proving that any positive $(a,abab\ldots)$-diagram over $\mathcal P$ can be composed from basic
diagrams together with an additional statement for a $(b,baba\ldots)$-diagram over $\mathcal P$. We prove
both facts simultaneously by induction on the number of cells in a diagram $\Phi$ with this property.

If $\Phi$ has no cells, there is nothing to prove. Let $\Phi$ have $a$ as a top label. Notice that the defining
relations of $\mathcal P$ always preserve the last letter of a word. So $a$ cannot be equal modulo this
presentation to a word that ends with $b$. Hence $\Phi$ is an $(a,(ab)^ma)$-diagram for some $m\ge1$. The top
cell of $\Phi$ has the form $a=(ba)^2$. Since the bottom path $p$ starts with $a$, the first letter $b$ of the
word $(ba)^2$ must correspond to the top path of a cell $b=(ab)^4$. These two cells form a basic diagram. So we
can cut it off. The rest will be a diagram with top path $a(ba)^5$ and bottom path $p$ labelled by $(ab)^ma$.

All vertices of a positive diagram belong to its bottom path. So it decomposes into a sum of diagrams for which
the top label of each of them is $a$ or $b$. If it is $a$, then the bottom label of a summand ends with $a$.
The length of the bottom path is odd so the bottom label has the form $(ab)^ka$ for some $k\ge0$. If the top
label of a summand is $b$, the same argument shows that the bottom label is of the form $(ba)^kb$. Thus all the
summands satisfy the inductive assumption (they have fewer cells than $\Phi$). Therefore they can be decomposed
into basic diagrams.

Now let $\Phi$ have $b$ as a top label. The top cell now is $b=(ab)^4$. The bottom path $p$ now starts with $b$.
Thus the first letter $a$ of $(ab)^4$ is the top path of a cell $a=(ba)^2$. The two cells together form a basic
diagram. We cut it off, and then repeat the same arguments as in the previous paragraph.

The image of the homomorphism is not Abelian. As above, we use the fact that generalized Thompson's groups $F_r$ have
no proper non-Abelian homomorphic images. Thus we have an isomorphism $F_{11}\cong G_{42}$.

The proof is complete.
\vspace{2ex}

Notice that the Fibonacci {\bf group} presented by $\mathcal P$ is a cyclic group $\mathbb Z_5$. The semigroup
with the same presentation is also finite, it has $10$ elements. However, for $n\ge5$ the Fibonacci semigroups presented by (\ref{fs})
turn out to be infinite. This makes unclear the structure of diagram groups $G_{n2}$ for that case (it is even possible that the groups
may be trivial). As for the generalization into another direction, we are able to describe completely the diagram groups over (\ref{pnr})
for the case $n=2$.

\begin{thm}
\label{johnn}
Let $s$ be a positive integer.

The diagram group with base $a$ over $\langle a,b\mid a=b(ab)^s,b=a(ba)^s\rangle$ is isomorphic to
generalized Thompson's group $F_{(2s+1)^2}$.

The diagram group with base $a$ over $\langle a,b\mid a=(ba)^s,b=(ab)^s\rangle$ is isomorphic to
generalized Thompson's group $F_{4s-1}$.

So the group $G_{2r}=\mathcal D(\mathcal P_{2r},a)$ is isomorphic to $F_{r^2}$ for odd $r$ and $F_{2r-1}$ for
even $r$, where $\mathcal P_{2r}=\langle a,b\mid a=ba\ldots,b=ab\ldots\rangle$ with the right-hand sides of the
defining relations of length $r\ge2$.
\end{thm}

\prf The case of odd $r=2s+1$ has the same proof as in Theorem~\ref{f9}. Basic diagrams here consist of $r+1$
cells. They correspond to the derivation $a=b(ab)^s$ with further replacements of all the $r$ letters of the
right-hand side according to the defining relators, and similarly for $b=a(ba)^s$ (we have a total symmetry here).
The bottom label of basic diagrams have length $r^2$. The proof goes without any changes for the general case.

Now let $r=2s$ be even. The construction of basic diagrams here is simpler. They consist of two cells only.
There is some similarity here to the construction from the proof of Theorem~\ref{f11}. Namely, we take the
cell $a=(ba)^s$ and replace the first letter in the right-hand side by $(ab)^s$. As a result, we get an
$(a,(ab)^{2s-1}a)$-diagram of two cells. We call it basic as well as the $(b,(ba)^{2s-1}b)$-diagram of two cells.
The bottom paths here have length $4s-1=2r-1$ so we are able to construct a homomorphism from $F_{2r-1}$ to the
diagram group and then show it is an isomorphism. The construction here is slightly easier than the one from
the proof of Theorem~\ref{f11} because of symmetry. This completes the proof.

\end{document}